\documentclass[a4paper]{article}
\usepackage{amsmath,amssymb,mathrsfs,graphicx,psfrag}
\DeclareMathOperator{\re}{Re}
\DeclareMathOperator{\im}{Im}
\begin{document}
\title{A method of fundamental solutions for doubly-periodic potential flow problems 
using the Weierstrass elliptic function}
\author{Hidenori Ogata%
\footnote{
Department of Computer and Network Engineering, 
Graduate School of Informatics and Engineering, 
The University of Electro-Communications, 
1-5-1 Chofugaoka, Chofu, Tokyo 182-8585, Japan. 
(e-mail) {\tt ogata@im.uec.ac.jp}
}}
\maketitle
\begin{abstract}
 In this paper, we propose a method of fundamental solutions for 
 the problems of two-dimensional potential flow past a doubly-periodic array of obstacles. 
 The solutions of these problems involve doubly-periodic functions, 
 and it is difficult to apply the conventional method of fundamental solutions to approximate them. 
 The method that we propose gives approximate solutions which is expressed by a linear combination of 
 periodic fundamental solutions constructed using the Weierstrass elliptic functions, 
 and it satisfies the periodicity that we expect. 
 Numerical examples show the effectiveness of our method. 
\end{abstract}

\medskip

\section{Introduction}
The method of fundamental solutions, or the charge simulation method 
\cite{FairweatherKarageorghis1998,Murashima1983}, is a fast solver of potential problems
\begin{equation*}
 \begin{cases}
  -\triangle u = 0 & \mbox{in} \ \mathscr{D} \\ 
  u = f & \mbox{on} \ \partial\mathscr{D}, 
 \end{cases}
\end{equation*}
where $\mathscr{D}$ is a domain in the $n$-dimensional Euclidean space $\mathbb{R}^n$ 
and $f$ is a function given on the boundary $\partial\mathscr{D}$. 
In the two-dimensional problems $(n=2)$, 
equalizing the Euclidean plane $\mathbb{R}^n$ with the complex plane $\mathbb{C}$, 
the method approximates the solution $u$ by 
\footnote{
The approximation shown here is that of the invariant scheme proposed by Murota 
\cite{Murota1993,Murota1995}. 
}
\begin{equation}
 \label{eq:mfs-solution0}
  u(z) \simeq u_N(z) = 
  Q_0 - \frac{1}{2\pi}\sum_{j=1}^{N}Q_j\log|z-\zeta_j|
  \quad ( \: z = x + \mathrm{i}y \: ),
\end{equation}
where $\zeta_1, \ldots, \zeta_N$ are points given in $\mathbb{C}\setminus\overline{\mathscr{D}}$ 
and $Q_0, Q_1, \ldots, Q_N$ are unknown real constants such that 
\begin{equation}
 \label{eq:q-sum-zero0}
  \sum_{j=1}^{N}Q_j = 0. 
\end{equation} 
We call $Q_j$ the ^^ ^^ charges'' and $\zeta_j$ the ^^ ^^ charge points''. 
We remark that the approximate solution $u_N$ exactly satisfies the Laplace equation in $\mathscr{D}$. 
Regarding the boundary condition, we pose the following collocation condition on $u_N$. 
\begin{equation}
 \label{eq:collocation-cond0}
  u_N(z_i) = f(z_i), \quad i = 1, \ldots, N,
\end{equation}
where $z_1, \ldots, z_N$ are points given on $\partial\mathscr{D}$. 
We call $z_i$ the ^^ ^^ collocation points''. 
The equation (\ref{eq:collocation-cond0}) is rewritten as
\begin{equation}
 \label{eq:collocation-cond02}
  Q_0 - \frac{1}{2\pi}\sum_{j=1}^{N}Q_j\log|z_i - \zeta_j| = f(z_i), 
  \quad i = 1, \ldots, N. 
\end{equation}
The equations (\ref{eq:q-sum-zero0}) and (\ref{eq:collocation-cond02}) form a system of 
linear equations for $Q_j$. 
We obtain $Q_j$ by solving the linear system (\ref{eq:q-sum-zero0}) and (\ref{eq:collocation-cond02}), 
and we obtain the approximate solution $u_N$. 
The method of fundamental solutions has the advantages that it is easy to program, its computational cost is low, 
and it achieves high accuracy such as exponential convergence under some conditions 
\cite{KatsuradaOkamoto1988,Katsurada1990,OgataKatsurada2014}. 
It was first used for electrostatic problems \cite{SingerSteinbiglerWeiss1974,Steinbigler-dissertation1969}, and now 
it is used in widely in science and engineering such as Helmholtz equation problems 
\cite{ChibaUshijima2009,OgataChibaUshijima2011} and studies on scattering of earthquake wave 
\cite{Sanchez-SezmaRosenblueth1979} and so on.

The method of fundamental solutions is also applied to the approximation of complex analytic functions. 
Let $f(z)$ be an complex analytic function in some complex domain $\mathscr{D}$ 
which is to be approximated and satisfies some boundary condition. 
Since the real part $\re f(z)$ is a harmonic function in $\mathscr{D}$, 
it can be approximated using the right hand side 
of (\ref{eq:mfs-solution0}), and the imaginary part $\im f(z)$, the conjugate harmonic function of $\re f(z)$, 
is approximated using 
\begin{equation*}
 - \frac{1}{2\pi}\sum_{j=1}^{N}Q_j\arg(z-\zeta_j). 
\end{equation*}
Then, the function $f(z)$ is approximated using the linear combination of the complex logarithmic functions 
\begin{equation}
 \label{eq:mfs-approx-analytic-function0}
 Q_0 - \frac{1}{2\pi}\sum_{j=1}^{N}Q_j\log(z - \zeta_j).
\end{equation}
From this point of view, Amano \cite{Amano1994,Amano1998} applied the method of fundamental solutions to 
numerical conformal mapping. 

\begin{figure}[htbp]
 \begin{center}
  \psfrag{D}{$\mathscr{D}$}
  \psfrag{w}{$\omega_1$}
  \psfrag{W}{$\omega_2$}
  \psfrag{0}{$D_{00}$}
  \psfrag{1}{$D_{10}$}
  \psfrag{2}{$D_{01}$}
  \psfrag{3}{$D_{11}$}
  \includegraphics[width=0.6\textwidth]{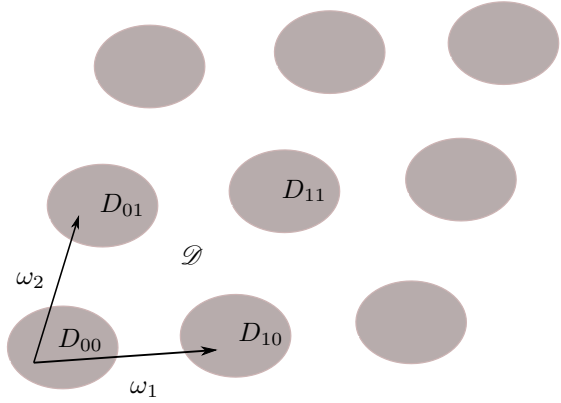}
 \end{center}
 \caption{Doubly-periodic array of obstacles.}
 \label{fig:periodic-domain}
\end{figure}

In this paper, we examine the problem of two-dimensional potential flow past a doubly-periodic array of obstacles 
as shown in Figure \ref{fig:periodic-domain}. 
A two-dimensional potential flow in a domain $\mathscr{D}$ is characterized by a complex velocity potential $f(z)$, 
a complex analytic function in $\mathscr{D}$ which gives the velocity field $\boldsymbol{v}=(u,v)$ by 
$f^{\prime}(z)=u-\mathrm{i}v$ and satisfies the boundary condition 
\begin{equation}
 \label{eq:boundary-cond-potential-flow0}
 \im f = \mbox{constant} \quad \mbox{on} \ \partial\mathscr{D}. 
\end{equation}
Physically, the condition (\ref{eq:boundary-cond-potential-flow0}) means that the fluid flows along the boundary 
$\partial\mathscr{D}$ since the contour lines of $\im f(z)$ give the streamlines. 
Therefore, we can obtain the complex velocity potential $f(z)$ by approximating it by the form 
(\ref{eq:mfs-approx-analytic-function0}). 
However, it is difficult to apply the method of fundamental solutions to our problem 
because the solution of our problem obviously involves a doubly-periodic function 
due to the periodicity of the problem, 
which it is difficult to approximate by the form (\ref{eq:mfs-approx-analytic-function0}) 
of the conventional method. 
To overcome this challenge, we propose a new method of fundamental solution method for our problem. 
In the proposed method, we approximate the solution involving a doubly-periodic function 
using a linear combination of periodic fundamental solutions, that is, complex logarithmic potentials with 
sources in a doubly-periodic array which is constructed by the Weierstrass sigma functions. 
Our method inherits the advantages of the conventional method of fundamental solutions 
and give an approximate solutions with the periodicity which we want. 
The author proposed a method of fundamental solution for the problem of two-dimensional potential flow with 
double periodicity, where we used the periodic fundamental solution constructed by the theta functions 
\cite{Ogata-arxiv2020-01}. 

The author has presented many works on methods of fundamental solution for periodic problems. 
In \cite{OgataOkanoAmano2002}, we presented a method of fundamental solutions for numerical conformal mappings 
of singly-periodic complex domains, where we approximated the mapping function using singly-periodic fundamental 
solutions, that is, the complex logarithmic potentials with sources in a singly-periodic array.  
The author proposed methods of fundamental solutions for the problems of Stokes flow past a periodic array of obstacles 
\cite{OgataAmanoSugiharaOkano2003,Ogata2006,OgataAmano2006,OgataAmano2010}, where the solutions are approximated using 
the periodic fundamental solutions of the Stokes equation, that is, the Stokes flow induced by concentrated forces 
in a periodic array. 
As to works on periodic Stokes flow, Zick and Homsy \cite{ZickHomsy1982} proposed an integral equation method 
for three-dimensional Stokes flow problems with a three-dimensional periodic array of spheres. 
Greengard and Kropinski \cite{GreengardKropinski2004} 
proposed an integral equation method for two-dimensional Stokes flow problems 
in doubly-periodic, where the approximate solution is given as a complex variable formulation 
and the fast multipole method is used. 
Liron \cite{Liron1978} proposed a study on Stokes flow due to infinite arrays of Stokeslets and its application to 
fluid transport by cilia. 
In addition, the author proposed a method of fundamental solutions for periodic plane elasticity \cite{Ogata2008}, 
where the solution is approximated using the periodic fundamental solutions of the elastostatic equation, that is, 
the displacements induced by concentrated forces in a periodic array. 

The remainder of this paper is structured as follows. 
Section \ref{sec:mfs} gives our method of fundamental solutions for periodic-potential flow. 
Section \ref{sec:example} presents some numerical examples which show the effectiveness of our method. 
In Section \ref{sec:conclusion}, we conclude this paper and presents some problems related to future studies. 
\section{Method of fundamental solutions}
\label{sec:mfs}
We examine the problem of two-dimensional potential flow past a doubly-periodic array of obstacles. 
In terms of mathematics, the flow domain is given by 
\begin{equation*}
 \mathscr{D} = \mathbb{C} \setminus \bigcup_{m,n\in\mathbb{Z}}\overline{D_{mn}}, 
\end{equation*}
where $D_{00}$ is one of obstacles, in terms of mathematics, a simply-connected domain in $\mathbb{C}$, and 
$D_{mn}$, $m,n\in\mathbb{Z}$ is given by 
\begin{equation*}
 D_{mn} = 
  \left\{ \: z + m\omega_1 + n\omega_2 \: | \: z \in D_{00} \: \right\}
\end{equation*}
with the complex number $\omega_1, \omega_2\in\mathbb{C}\setminus\{0\}$ such that 
$\im(\omega_2/\omega_1)>0$ giving the periods of our problem. 
The complex velocity potential of our problem is a complex analytic function $f(z)$ in the domain 
$\mathscr{D}$ satisfying the boundary condition
\begin{equation}
 \label{eq:boundary-cond}
  \im f(z) = \mbox{constant} \quad \mbox{on} \ \partial D_{mn}, \ m,n \in\mathbb{Z}. 
\end{equation}
Additionally, we assume that the average of the flow is the uniform flow, in whose direction 
we take the real axis, when it is observed along the direction of $\omega_1$, that is, 
\begin{equation}
 \label{eq:average-flow-cond0}
  \int_{C}\boldsymbol{v}\cdot\mathrm{d}\boldsymbol{r} = U(\re\omega_1), 
  \quad 
  \int_{C}\boldsymbol{v}\cdot\boldsymbol{n}\mathrm{d}s = U(\im\omega_1), 
\end{equation}
where $U$ is the magnitude of the uniform flow, 
$C$ is the curve in $\mathscr{D}$ connecting two points $z_0, z_0+\omega_1\in\mathscr{D}$
and $\boldsymbol{n}$ is the unit normal vector on $C$. 
Using $f^{\prime}(z) = u - \mathrm{i}v$, the condition (\ref{eq:average-flow-cond0}) is rewritten as 
\begin{equation}
 \label{eq:average-flow-cond}
 \int_{z_0}^{z_0+\omega_1}\mathrm{d}f(z) = 
 f(z_0 + \omega_1) - f(z_0) = 
 U\omega_1.
\end{equation}
Therefore, our problem is to find a complex velocity potential $f(z)$, 
an analytic function in $\mathscr{D}$ such that it satisfies the boundary condition 
(\ref{eq:boundary-cond}) and the condition (\ref{eq:average-flow-cond}). 

We propose to approximate the complex velocity potential $f(z)$ in the following form 
according to \cite{Hasimoto2008}.
\begin{equation}
 \label{eq:periodic-mfs-solution}
 f_N(z) = 
 \left( U - \frac{\mathrm{i}\eta_1}{2\pi\omega_1}\sum_{j=1}^{N}Q_j \zeta_j\right)z 
 - 
 \frac{\mathrm{i}}{2\pi}\sum_{j=1}^{N}Q_j\log\sigma(z - \zeta_j),
\end{equation}
where $\zeta_1, \ldots, \zeta_N$ are points given in $D_{00}$, $Q_1, \ldots, Q_N$ are unknown real coefficients 
such that 
\begin{equation}
 \label{eq:q-sum-zero}
  \sum_{j=1}^{N}Q_j = 0,
\end{equation}
$\sigma(z)$ is the Weierstrass sigma function \cite{Armitage-Eberlein2006}, 
and $\eta_1$ is given by 
$\eta_1 = \zeta(\omega_1)$ using the Weierstrass zeta function $\zeta(z)$. 
We call the real coefficients $Q_j$ the ^^ ^^ charges'' and $\zeta_j$ the ^^ ^^ charge points''. 
Since $\sigma(z)$ is an entire function with simple zeros at $z=m\omega_1+n\omega_2$, $m,n\in\mathbb{Z}$, 
the functions $\sigma(z-\zeta_j)$ appearing on the right hand side of (\ref{eq:periodic-mfs-solution}) 
are complex logarithmic potential with sources in a doubly-periodic array. 
The approximate potential $f_N(z)$ satisfies the condition (\ref{eq:average-flow-cond}). 
In fact, we have
\begin{align*}
 & f_N(z_0 + \omega_1) - f_N(z_0) 
 \\
 = \: & 
 \left( U - \frac{\mathrm{i}\eta_1}{2\pi\omega_1}\sum_{j=1}^{N}Q_j\zeta_j\right)\omega_1 
 \\
 & 
 - \frac{\mathrm{i}}{2\pi}\sum_{j=1}^{N}Q_j
 \left\{
 \log\left(-\exp\left[\eta_1\left(z-\zeta_j+\frac{\omega_1}{2}\right)\right]\log\sigma(z-\zeta_j)\right)
 - 
 \log\sigma(z-\zeta_j)
 \right\}
 \\ 
 = \: & 
 \left( U - \frac{\mathrm{i}\eta_1}{2\pi\omega_1}\sum_{j=1}^{N}Q_j\zeta_j\right)\omega_1 
 - \frac{\mathrm{i}}{2\pi}\sum_{j=1}^{N}Q_j
 \left\{ \log(-1) + \eta_1\left(z - \zeta_j + \frac{\omega_1}{2}\right)\right\}
 \\
 = \: & 
 U\omega_1, 
\end{align*}
where we used the pseudo-periodicity of the sigma function 
\begin{equation}
 \label{eq:pseudo-periodicity-sigma}
  \sigma(z + \omega_k) = 
  - \exp\left[\eta_k\left(z + \frac{\omega_k}{2}\right)\right]\sigma(z), \quad k = 1, 2
\end{equation}
with $\eta_k = \zeta(\omega_k)$, $k=1,2$ and (\ref{eq:q-sum-zero}). 
The approximate potential $f_N(z)$ satisfies the pseudo-periodicity
\begin{equation}
 \label{eq:pseudo-periodicity-mfs-solution}
  f_N(z+\omega_1) = f_N(z) + U\omega_1, \quad 
  f_N(z+\omega_2) = f_N(z) + U\omega_2 + \frac{1}{\omega_1}\sum_{j=1}^{N}Q_j\zeta_j
\end{equation}
due to (\ref{eq:pseudo-periodicity-sigma}) and (\ref{eq:q-sum-zero}). 
Then, the complex velocity $f_N^{\prime}(z) = u_N(z) - \mathrm{i}v_N(z)$ satisfies 
\begin{equation*}
 f_N^{\prime}(z + \omega_k) = f_N^{\prime}(z), \quad k = 1, 2, 
\end{equation*}
that is, it is an elliptic function with periods $\omega_1$ and $\omega_2$. 
Regarding the boundary condition (\ref{eq:boundary-cond}), 
we pose the following collocation condition on $f_N(z)$
\begin{equation}
 \label{eq:collocation-cond}
  \im f_N(z_i) = C, \quad i = 1, \ldots, N,
\end{equation}
where $z_1, \ldots, z_N$ are points given on $\partial D_{00}$, 
and $C$ is a real constant. 
We call $z_i$ the ^^ ^^ collocation condition''. 
The collocation condition (\ref{eq:collocation-cond}) is rewritten as
\begin{multline}
 \label{eq:collocation-cond2}
 - C - \frac{1}{2\pi}\sum_{j=1}^{N}Q_j
 \left\{ \log|\sigma(z_i - \zeta_j)| + \frac{\eta_1}{\omega_1}\re\left(\sum_{j=1}^{N}\zeta_jz_i\right) \right\}
 = - U\im z_i, 
 \\ 
 i = 1, \ldots, N, 
\end{multline}
which form a system of linear equations for $C$ and $Q_j$ together with (\ref{eq:q-sum-zero}). 
We obtain the charges $Q_j$ by solving the linear system (\ref{eq:q-sum-zero}) and (\ref{eq:collocation-cond2}) 
and obtain the approximate potential $f_N(z)$. 
By (\ref{eq:collocation-cond}), $f_N(z)$ approximately satisfies the boundary condition (\ref{eq:boundary-cond}) 
on the boundaries of other obstacles $D_{mn}$, $m,n\in\mathbb{Z}$ due to the pseudo-periodicity 
(\ref{eq:pseudo-periodicity-mfs-solution}). 
\section{Numerical examples}
\label{sec:example}
In this section, we show some numerical examples which show the effectiveness of our method. 
All the computations were performed using programs coded in C++ with double precision. 

\begin{figure}[htbp]
 \begin{center}
  \begin{tabular}{cc}
   \psfrag{x}{$\re z/r$}
   \psfrag{y}{$\im z/r$}
   \includegraphics[width=0.49\textwidth]{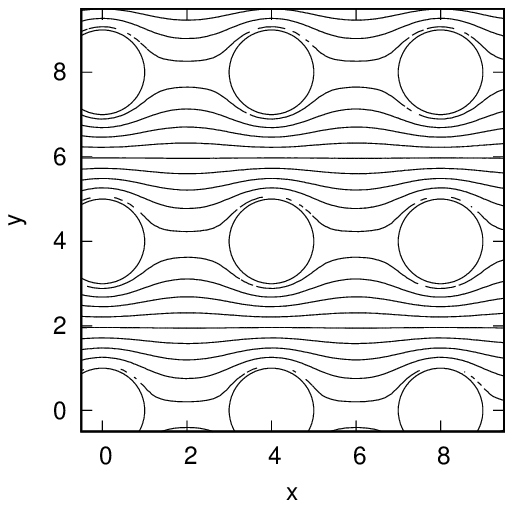} & 
       \psfrag{x}{$\re z/r$}
       \psfrag{y}{$\im z/r$}
       \includegraphics[width=0.49\textwidth]{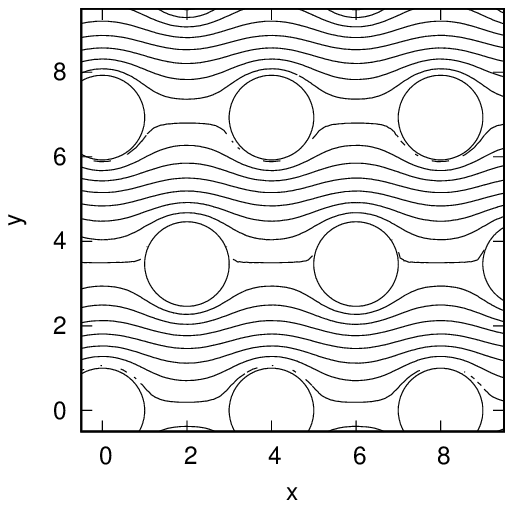}
       \\
   $(\omega_1, \omega_2) = (4r,4r\mathrm{i})$ & 
       $(\omega_1, \omega_2) = (4r,4r\mathrm{e}^{\mathrm{i}\pi/3})$
       \\
   \psfrag{x}{$\re z/r$}
   \psfrag{y}{$\im z/r$}
   \includegraphics[width=0.49\textwidth]{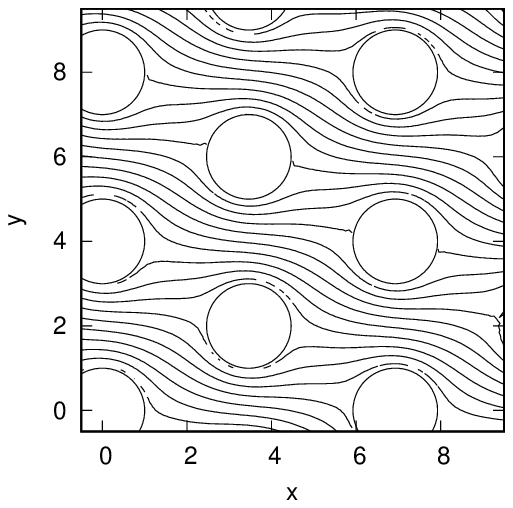} & 
       \psfrag{x}{$\re z/r$}
       \psfrag{y}{$\im z/r$}
       \includegraphics[width=0.49\textwidth]{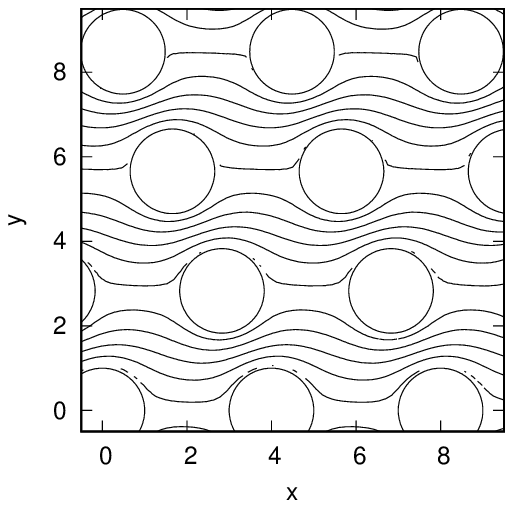}
       \\ 
   $(\omega_1,\omega_2) = (4r\mathrm{e}^{\mathrm{i}\pi/6}, 4r\mathrm{i})$ & 
       $(\omega_1,\omega_2) = (4r, 4r\mathrm{e}^{\mathrm{i}\pi/4})$
  \end{tabular}
 \end{center}
 \caption{The streamlines of two-dimensional potential flows 
 past a doubly-periodic array of cylinders. }
 \label{fig:example1}
\end{figure}

We computed the two-dimensional potential flow past a doubly-periodic array of cylinders 
of radius $r>0$
\begin{gather*}
 \mathscr{D} = \mathbb{C} \setminus \bigcup_{m,n\in\mathbb{Z}}\overline{D_{mn}}, 
 \intertext{where}
 D_{mn} = 
 \left\{ \: z \in\mathbb{C} \: | \: |z - m\omega_1 - n\omega_2| < r 
 \: \right\}, 
 \quad m, n \in\mathbb{Z}, 
\end{gather*}
and $\omega_1, \omega_2$ are the periods of the array such that 
\begin{math}
 \im(\omega_2/\omega_1) > 0,  
\end{math}
by our method. 
In our method, we took the charge points $\zeta_j$ and the collocation points $z_i$ respectively as
\begin{equation}
 \label{eq:charge-collocation-point}
  \zeta_j = qr\exp\left(\mathrm{i}\frac{2\pi(j-1)}{N}\right), \quad 
  z_j = r\exp\left(\mathrm{i}\frac{2\pi(j-1)}{N}\right), \quad j = 1, \ldots, N, 
\end{equation}
where $q$ is a constant such that $0 < q < 1$, which is taken as $q=0.7$ in the example. 
Figure \ref{fig:example1} shows the streamline of the flows for 
some pairs of $(\omega_1, \omega_2)$

\begin{figure}[htbp]
 \begin{center}
  \begin{tabular}{cc}
   \psfrag{n}{$N$}
   \psfrag{e}{$\epsilon_N$}
   \includegraphics[width=0.49\textwidth]{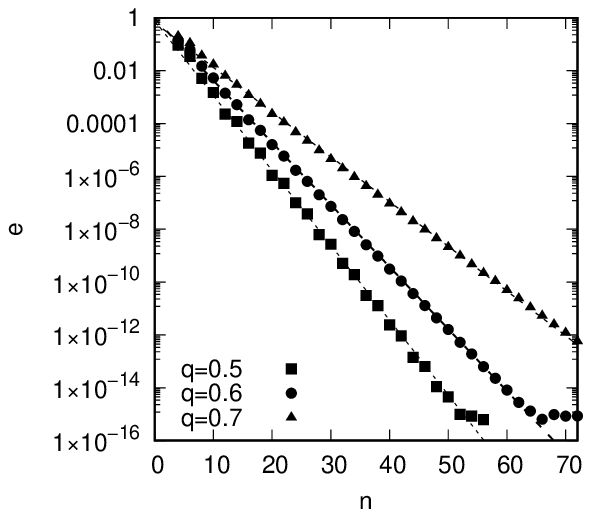} & 
       \psfrag{n}{$N$}
       \psfrag{e}{$\epsilon_N$}
       \includegraphics[width=0.49\textwidth]{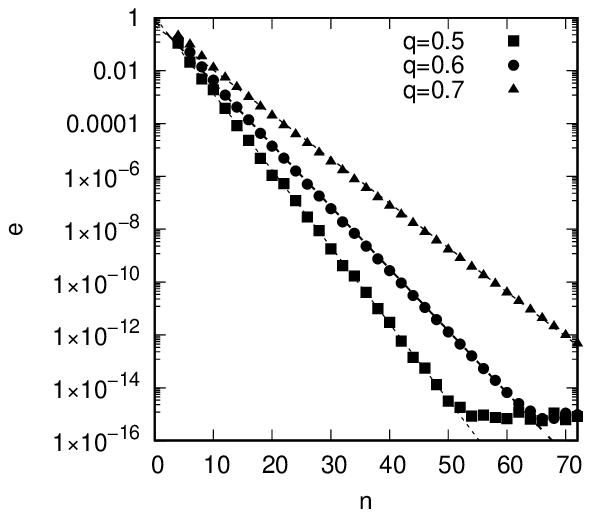}
       \\
   $(\omega_1, \omega_2) = (4r,4r\mathrm{i})$ & 
       $(\omega_1, \omega_2) = (4r, 4r\mathrm{e}^{\mathrm{i}\pi/3})$
       \\
   \psfrag{n}{$N$}
   \psfrag{e}{$\epsilon_N$}
   \includegraphics[width=0.49\textwidth]{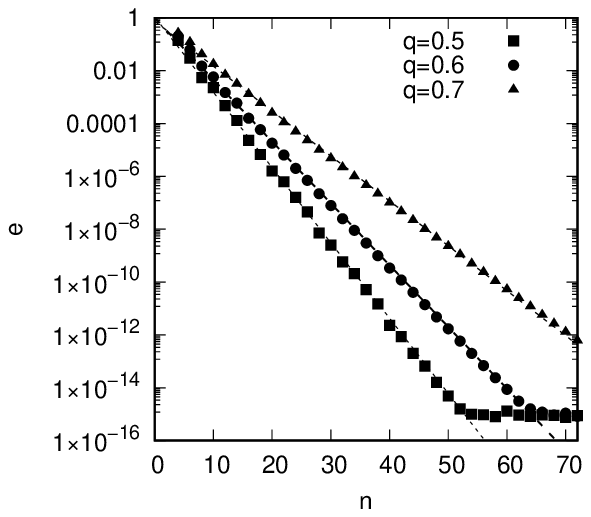} &
       \psfrag{n}{$N$}
       \psfrag{e}{$\epsilon_N$}
       \includegraphics[width=0.49\textwidth]{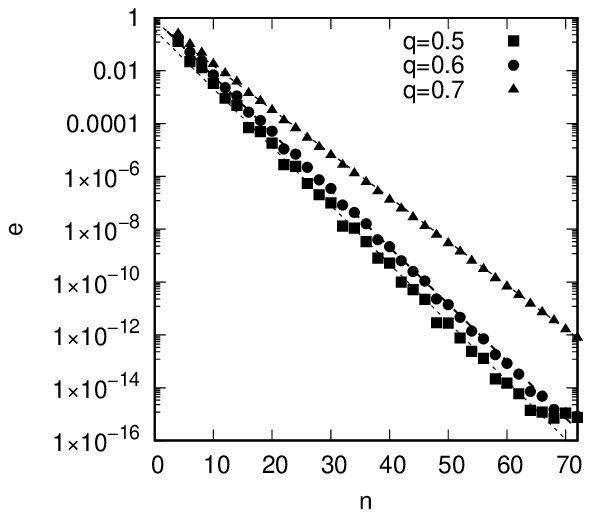}
       \\
   $(\omega_1,\omega_2)=(4r\mathrm{e}^{\mathrm{i}\pi/6},4r\mathrm{i})$ & 
       $(\omega_1,\omega_2)=(4r,4r\mathrm{e}^{\mathrm{i}\pi/4})$
  \end{tabular}
 \end{center}
 \caption{The error estimate $\epsilon_N$ of our method. }
 \label{fig:error1}
\end{figure}

To estimate the accuracy of our method, we computed 
\begin{equation*}
 \epsilon_N = 
  \frac{1}{Ur}\max_{z\in\partial D_{00}}|\im f_N(z) - C|, 
\end{equation*}
where $C$ is the constant appearing in (\ref{eq:collocation-cond2}). 
The value $\epsilon_N$ shows how accurately the approximate potential $f_N(z)$ 
satisfies the boundary condition (\ref{eq:boundary-cond}). 
Figure \ref{fig:error1} shows the value $\epsilon_N$ computed for the example with the parameter $q$ 
in (\ref{eq:charge-collocation-point}) taken as some values. 
The figures show that $\epsilon_N$ decays exponentially as the number of unknowns $N$ increases. 
Table \ref{tab:error1} shows the decay rates of $\epsilon_N$ computed by the least square fitting 
using the {\tt fit} command of the software {\it gnuplot}. 
The table shows that the decay rate of $\epsilon_N$ roughly obeys the rule 
\begin{equation*}
 \epsilon_N = \mathrm{O}(q^N).
\end{equation*}

\begin{table}[htbp]
 \caption{The behavior of the error estimate $\epsilon_N$.}
 \begin{center}
  \begin{tabular}{cccc}
   \hline \rule{0pt}{12pt} 
   $(\omega_1, \omega_2)$ & 
   \multicolumn{3}{c}{$q$}
       \\
   \cline{2-4} \rule{0pt}{12pt}
   & $0.5$ & $0.6$ & $0.7$ \\
   \hline \rule{0pt}{12pt}
   $(4r,4r\mathrm{i})$ & $\mathrm{O}(0.52^N)$ & $\mathrm{O}(0.58^N)$ & $\mathrm{O}(0.67^N)$ \\ 
   $(4r,4r\mathrm{e}^{\mathrm{i}\pi/3})$ & $\mathrm{O}(0.51^N)$ & $\mathrm{O}(0.58^N)$ & $\mathrm{O}(0.68^N)$ \\ 
   $(4r\mathrm{e}^{\mathrm{i}\pi/6},4r\mathrm{i})$ & $\mathrm{O}(0.52^N)$ & $\mathrm{O}(0.58^N)$ & $\mathrm{O}(0.68^N)$
   \\
   $(4r,4r\mathrm{e}^{\mathrm{i}\pi/4})$ & $\mathrm{O}(0.60^N)$ & $\mathrm{O}(0.61^N)$ & $\mathrm{O}(0.68^N)$ \\ 
   \hline 
  \end{tabular}
 \end{center}
 \label{tab:error1}
\end{table}
\section{Concluding Remarks}
\label{sec:conclusion}
In this paper, we examined the problems of two-dimensional potential flow past a doubly-periodic array of 
obstacles and proposed a method of fundamental solutions for these problems. 
It is difficult to apply the conventional method of fundamental solution to our problems 
because the solution involves a doubly-periodic functions. 
We proposed a method of fundamental solutions for our problems, where the solution is approximated 
using the periodic fundamental solutions, that is, the complex logarithmic 
potentials with sources in a doubly-periodic array and constructed by the Weierstrass sigma functions. 
The proposed method inherits the advantages of the conventional method and approximates well the solution 
involving a periodic function. 
The numerical examples showed the effectiveness of our method. 

We have two issues regarding this paper for future studies. 
The first problem is to extend our method to other periodic problems than the Laplace equation such as 
the Stokes equation. 
In the author's previous works on periodic Stokes flow 
\cite{OgataAmanoSugiharaOkano2003,Ogata2006,OgataAmano2006}, the solutions are approximated using the periodic 
fundamental solutions which was presented by Hasimoto \cite{Hasimoto1959} and expressed by a Fourier series. 
In these years, Hasimoto presented the periodic fundamental solutions of the Stokes equation using 
the Weierstrass elliptic functions \cite{Hasimoto2008}, 
and it is interesting to construct a method of fundamental 
solutions using these periodic fundamental solutions.

The second problem is a theoretical study on the accuracy of our method. 
Theoretical error estimates on method of fundamental solutions are presented for special problems such as 
two-dimensional potential or Helmholtz equation problems in a disk 
\cite{KatsuradaOkamoto1988,ChibaUshijima2009,OgataChibaUshijima2011} and two-dimensional potential problems 
in a domain with an analytic boundary \cite{Katsurada1990,OgataKatsurada2014}; however, 
there still remain many problems including the problem of this paper to which theoretical error estimates 
are to be given. 
This is one of the most important problems on the method of fundamental solutions. 
\bibliographystyle{plain}
\bibliography{arxiv2020_02}
\end{document}